\documentclass[12pt]{article}

\usepackage[margin=2.5cm, top=2.5cm, bottom=2.5cm]{geometry}

\usepackage[utf8]{inputenc} 
\usepackage[T1]{fontenc}    
\usepackage{url}            
\usepackage{booktabs}       
\usepackage{amsfonts}       
\usepackage{nicefrac}       
\usepackage{microtype}      
\usepackage[authoryear]{natbib}
\usepackage{amsmath,amssymb,amsthm,xcolor}
\usepackage{array}
\usepackage{float}
\usepackage[hidelinks,hypertexnames=false]{hyperref}
\usepackage[capitalize]{cleveref}
\theoremstyle{definition}
\newtheorem{example}{Example}[section]

\theoremstyle{plain}
\newtheorem{theorem}[example]{Theorem}

\newtheorem{proposition}[example]{Proposition}

\usepackage{graphicx}
\usepackage{subcaption}
\usepackage{mdframed}
\usepackage{lipsum}
\usepackage{wrapfig}
\usepackage{tikz}
\usetikzlibrary{arrows.meta}
\usepackage{enumitem}
\allowdisplaybreaks
\usepackage{epstopdf}
\usepackage{booktabs}
\usepackage{changepage}
\usepackage{multirow}
\usepackage{algorithm}
\usepackage{algorithmic}
\usepackage{mathtools}
\crefname{assumption}{Assumption}{Assumptions}

\numberwithin{equation}{section}
\usepackage[nottoc]{tocbibind}
\usepackage{mathrsfs}
\usepackage{stmaryrd}

\usepackage{amsmath,amsfonts,bm}

\def\eqref#1{equation~\ref{#1}}

\def\1{\bm{1}}

\def\vzero{{\bm{0}}}

\def\va{{\bm{a}}}
\def\vb{{\bm{b}}}

\def\vu{{\bm{u}}}
\def\vv{{\bm{v}}}

\def\vx{{\bm{x}}}
\def\vy{{\bm{y}}}

\def\eva{{a}}
\def\evb{{b}}

\def\evu{{u}}
\def\evv{{v}}

\def\evx{{x}}

\def\mA{{\bm{A}}}

\def\mH{{\bm{H}}}
\def\mI{{\bm{I}}}

\def\mM{{\bm{M}}}

\def\mO{{\bm{O}}}

\DeclareMathAlphabet{\mathsfit}{\encodingdefault}{\sfdefault}{m}{sl}
\SetMathAlphabet{\mathsfit}{bold}{\encodingdefault}{\sfdefault}{bx}{n}
\newcommand{\tens}[1]{\bm{\mathsfit{#1}}}

\def\tS{{\tens{S}}}
\def\tT{{\tens{T}}}
\def\tU{{\tens{U}}}

\def\tones{{\tens{1}}}

\newcommand{\Var}{\mathrm{Var}}

\title{\LARGE Global convergence of gradient descent for phase retrieval}

\begin{document}

\author{\large Théodore Fougereux\thanks{\url{theodore.fougereux@polytechnique.edu}, Centre de Mathématiques Appliquées, Ecole Polytechnique, Palaiseau, France.} \and Cédric Josz\thanks{\url{cj2638@columbia.edu}, IEOR, Columbia University, New York. Research supported by NSF EPCN grant 2023032 and ONR grant N00014-21-1-2282.} \and Xiaopeng Li\thanks{\url{xl3040@columbia.edu}, IEOR, Columbia University, New York.}}
\date{}

\maketitle
\vspace*{-5mm}
\begin{center}
    \textbf{Abstract}
    \end{center}
    \vspace*{-4mm}
 \begin{adjustwidth}{0.2in}{0.2in}

We propose a tensor-based criterion for benign landscape in phase retrieval and establish boundedness of gradient trajectories. This implies that gradient descent will converge to a global minimum for almost every initial point.


\end{adjustwidth} 
\vspace*{3mm}
\vspace*{-3mm}

\section{Introduction}
We study the problem of phase retrieval, which aims to recover a complex valued vector $\vx^{\natural}\in\mathbb{C}^n$ from its intensity measurements 
\begin{equation}\label{eq:intensity_measure}
    y_i=|\langle \va_i, \vx^{\natural}\rangle|^2, \quad i=1,\ldots,m,
\end{equation}
where $\va_i\in\mathbb{C}^n$, $i=1,\ldots,m$, are known complex vectors and $m$ is the number of measurements. This problem has attracted high interest due to its broad applications in X-ray crystallography \citep{elser2018benchmark}, microscopy \citep{miao2008extending}, astronomy \citep{fienup1987phase} and optical imaging \citep{shechtman2015phase}. 

The phase retrieval problem is NP-hard if only very few measurements, e.g. $m=n+1$ \citep{fickus2014phase}, are given. However, a wide range of algorithms can recover $\vx^{\natural}$ up to a global phase shift provided enough measurements. Early methods with provable performance guarantees usually formulate it into a convex optimization problem, such as a semidefinite programming problem \citep{candes2013phaselift,candes2015phaseR,waldspurger2015phase,mcrae2022optimal} or basis pursuit problem \citep{goldstein2018phasemax}. These methods are usually computationally challenging in high-dimensional cases. To address this issue, more recent works take nonconvex approaches, such as alternating minimization \citep{wen2012alternating,netrapalli2013phase,waldspurger2018phase,zhang2020phase}; gradient descent type algorithms, including Wirtinger flow \citep{candes2015phase,chen2015solving,ma2020implicit}, truncated amplitude flow \citep{wang2017solving}, vanilla gradient descent \citep{chen2019gradient}, Riemannian gradient descent \citep{cai2024solving}; and Newton type algorithms \citep{gao2017phaseless,ma2018globally}. 

Convex methods mentioned above can achieve optimal sample complexity $m=O(n)$, but require an initialization close enough to the ground truth $\vx^{\natural}$. For nonconvex algorithms, $O(n)$ sample complexity can be achieved with a careful initialization \citep{chen2015solving,wang2017solving,waldspurger2018phase,cai2024solving}. However, such initialization can be computationally inefficient when the dimension $n$ is large. In practice, random initialization is more plausible. In the random initialization regime, it is known that $O(n\log^{13} n)$ samples are sufficient to guarantee a nearly linear convergence rate for vanilla gradient descent \citep{chen2019gradient}. To obtain a lower sample complexity in this regime, a benign global landscape result was analyzed for both the intensity measurement \cref{eq:intensity_measure} \citep{sun2018geometric,cai2023nearly} and the amplitude measurement $y_i=|\langle \va_i, \vx^{\natural}\rangle|$ \citep{cai2022solving} for $\ell_2$-loss function. Both the intensity and the amplitude model have no spurious local minima and saddle points are strict, provided $O(n\log n)$ and $O(n)$ samples respectively. For intensity model with $o(n\log n)$ samples, it is suggested in \cite{liu2023local} that the benign landscape is hard to obtain due to the lack of convexity around the ground truth. 
Nevertheless, no global convergence guarantee of gradient descent algorithm can be easily deduced because the Lipschitz gradient assumption in the classical convergence results \citep{lee2016gradient} does not hold, and the bounded iterates assumptions are hard to verify.

In this paper, we aim at showing the global convergence of vanilla gradient descent algorithm with arbitrary initialization for general phase retrieval problem. We propose a new tensor based criterion to guarantee a benign global landscape for phase retrieval problem. Furthermore, we show that such objective function has bounded gradient trajectories, which give us the certificate to apply the general global convergence result in \citep{josz2023global}. 

The previous version of this manuscript proposed a sample complexity which relied on a tensor concentration inequality \cite[Theorem 3]{even2021concentration} which is false (see the arxiv version for a corrected and weaker inequality). In particular, the sample complexity would contradict the landscape result of \cite[Corollary 2.2]{liu2023local}. \\\\

\noindent\textbf{Acknowledgements} We wish to thank Kaizhao Liu and 
Andrew Duncan Mc Rae for their valuable feedback on the first version of this manuscript. 




\section{Problem formulation and related work}

We first present the problem formulation for phase retrieval and the algorithm to solve it, and then review related work and provide a comparative overview of state-of-the-art convergence results and our new result. 

\textbf{Notation:} Denote $\mathbb{R}$ and $\mathbb{C}$ as the real and complex fields respectively. Denote $\|\vx\|$ as the $\ell_2$-norm of vector $\vx\in\mathbb{C}^n$. For $\va,\vb\in\mathbb{C}^n$, denote $\langle \va,\vb\rangle=\sum_{i=1}^n\eva_i\overline{\evb_i}$, where $\overline{\evb_i}$ is the complex conjugate of $\evb_i$. We denote $\va\otimes \vb$ as the tensor product of $\va$ and $\vb$ and $\va^{\otimes r}$ as the tensor product of $r$ vectors $\va$. 
Denote $\mI_n$ as the $n\times n$ identity matrix and $\mO_n$ as the $n\times n$ zero matrix. We call $\tT=\vu_1\otimes\cdots\otimes \vu_p$ a tensor of rank 1 and order $p$. Denote $\|\tT\|=\|\vu_1\|\cdots\|\vu_p\|$ and $\|\tT\|_{\rm op}=\sup_{\|\vx_1\otimes\cdots\otimes\vx_p\|=1}\langle \tT, \vx_1\otimes\cdots\otimes\vx_p\rangle$ as its operator norm. Denote $\mH_f(\vx)$ the Hessian of a function $f$ at a point $\vx$.  We write $f(n)=O(g(n))$ if $f(n)\leq C_1g(n)$, $f(n)=\Omega(g(n))$ if $f(n)\geq C_2g(n)$, and $f(n)=\Theta(g(n))$ if $C_2g(n)\leq f(n)\leq C_1g(n)$ for some constants $C_1,C_2>0$. 

\subsection{Problem formulation}
We denote $\vx^{\natural}\in\mathbb{C}^n$ as the ground truth vector we want to recover. Intensity-based phase retrieval problem aims at minimizing the empirical $\ell_2$-loss of $m$ intensity measurements $y_i=|\langle \va_i, \vx^{\natural}\rangle|^2$, 
\begin{align}\label{eq:in_phase_retrieval}
    \min_{\vx\in\mathbb{C}^n}\quad f(\vx):=\sum_{i=1}^{m} \left(|\langle \va_i,\vx\rangle|^2-y_i\right)^2, 
\end{align}
Since all vectors of the form $e^{i\theta}\vx^{\natural}$ are  global minimizers of $f$, we can only expect to recover $\vx^{\natural}$ up to a phase shift by solving \cref{eq:in_phase_retrieval}. For the purpose of comparison, we also recall the amplitude-based formulation mentioned in the introduction, 
\begin{equation*}
    \min_{\vx\in\mathbb{C}^n}\quad \sum_{i=1}^{m} \left(|\langle \va_i,\vx\rangle|-\sqrt{y_i}\right)^2. 
\end{equation*}
It is easy to see that the intensity-based formulation has a smooth objective function while the amplitude-based formulation has a nonsmooth one.

In order to apply vanilla gradient descent to solve \cref{eq:in_phase_retrieval} as it is applied in the real case \citep{chen2019gradient}, we introduce the following mapping for any $\vv\in\mathbb{C}^n$: 
\begin{equation}\label{eq:mapping}
    \vv \mapsto \vv^+ := \begin{bmatrix}
        {\rm Re}(\vv) \\ {\rm Im}(\vv)
    \end{bmatrix}, \quad \vv \mapsto \vv^- := \mM\vv^+ = \begin{bmatrix}
        -{\rm Im}(\vv) \\ {\rm Re}(\vv)
    \end{bmatrix} \text{ where } \mM:=\begin{bmatrix}
        \mO_n & -\mI_n \\
        \mI_n & \mO_n
    \end{bmatrix}.
\end{equation}
By using \cref{eq:mapping}, we can define $\va_i^+$, $\va_i^-$, $\vx^+$, $\vx^-$, $\vx^{\natural+}$ and $\vx^{\natural-}$. Then, with a little abuse of notation over $f$, we can obtain a equivalent form of $f$ in terms of $\vx^+$ as 
\begin{equation}\label{eq:realPR}
    f(\vx) = \sum_{i=1}^{m} \left(|\langle \va_i,\vx\rangle|^2-y_i\right)^2 = \sum_{i=1}^{m} \left(\langle \mA_i\vx^+,\vx^+\rangle-y_i\right)^2 := f(\vx^+),
\end{equation}
where $\mA_i=\va_i^+(\va_i^+)^{\rm T}+\va_i^-(\va_i^-)^{\rm T}$. Now we can treat $f$ as a function from $\mathbb{R}^{2n}$ to $\mathbb{R}$ w.r.t. $\vx^+$. 

The gradient descent algorithm to minimize $f$ is given by 
\begin{equation}\label{eq:gd_realPR}
    \vx^+_{k+1} = \vx^+_k - \alpha_k \nabla_{\vx^+} f(\vx^+_k), \quad \forall\, k\in\mathbb{N},
\end{equation}
where $\vx^+_0\in\mathbb{R}^{2n}$ is any given initial point and $\alpha_k>0$ is any step size. The goal of this paper is to study the global convergence and sample complexity of gradient descent algorithm \cref{eq:gd_realPR} for $f$ defined in \cref{eq:realPR}. As we did in \cref{eq:gd_realPR}, all gradient and Hessian in the rest of this paper will be taken w.r.t. $\vx^+$.

\subsection{Related work}

Global benign landscape results usually refer to properties that all local minima of a function are global minima and all saddle points are strict, i.e., Hessian matrix having a negative eigenvalue at the saddle point. For phase retrieval problem, existing results obtain $O(n)$ sample complexity when the objective function is quadratic-like far from the origin, including amplitude model (piecewise quadratic) \citep{cai2022solving}, quotient intensity model (quartic divided by quadratic) \citep{cai2022global}, and perturbed amplitude model (truncated quartic) \citep{cai2021global}. For purely quartic objective function, like intensity model, people usually require $O(n\,\text{\rm poly}(\log n))$ samples to obtain the global benign landscape result \citep{sun2018geometric,cai2023nearly}.  \Cref{table:landscape} summarizes the above results. 

\begin{table}[H]
\centering
\begin{tabular}{@{}lccc@{}}
\hline
\textbf{Work} & \textbf{Objective function} & \textbf{Sample complexity} & \textbf{Probability}\\ \hline
\citet{geometricanalisys2018} & quartic & $O(n \log^3 n)$ & $1-\Theta(1/m)$\\
\citet{cai2021global} & truncated quartic & $O(n)$ & $1-O(1/m^2)$ \\
\cite{cai2022global} & quartic over quadratic & $O(n)$ & $1-\exp(-\Theta(m))$ \\
\cite{cai2022solving} & piecewise quadratic & $O(n)$ & $1-\exp(-\Theta(m))$ \\ 
\cite{cai2023nearly} & quartic & $O(n\log n)$ & $1-\Theta(1/m)$ \\ \hline
\end{tabular}
\caption{Global benign landscape results of different objective functions and sample complexity in literature. Here ``quartic'' and ``quadratic'' means 4th order and 2nd order polynomial respectively. }\label{table:landscape}
\end{table}

Convergence results for gradient descent type algorithms to solve phase retrieval problem either requires special initialization that are usually already close a ground truth \citep{candes2015phase,chen2015solving,wang2017solving,cai2024solving} or suboptimal sample complexity $O(n\log^{13} n)$ \citep{chen2019gradient}, as summarized in \Cref{table:convergence}. In the former case, the iterates start nearby a ground truth, where a desirable local landscape occurs and is enough to guarantee a contraction property of the distance of iterates towards ground truth. In the latter case, a careful analysis of some nonlinear approximate dynamics and complicated concentration results are required to show that after a short period of time, the iterates will enter the neighborhood with  contraction property. Contrary to the existing proof idea, we develop a new and concise proof of global convergence results for phase retrieval problem by connecting global benign landscape results and continuous time gradient dynamics results. 

\begin{table}[!ht]
\centering
\begin{tabular}{@{}lccc@{}}
\hline
\textbf{Work} & \textbf{Algorithm} & \textbf{Initialization} &  \textbf{Sample complexity}  \\ \hline
\citet{candes2015phase} & WF & spectral & $O(n\log n)$ \\
\cite{chen2015solving} & TWF & spectral & $O(n)$ \\
\cite{wang2017solving} & TAF & orthogonality-promoting & $O(n)$ \\
\cite{chen2019gradient} & GD & $\mathcal{N}(\vzero,n^{-1}\mI_n)$ & $O(n\log^{13} n)$\\
\cite{cai2024solving} & TRGrad & close to ground truth & $O(n)$ \\ \hline
\end{tabular}
\caption{Global convergence results of gradient descent type algorithms with different initialization and sample complexity in literature. Here ``GD'' stands for vanilla gradient descent, ``WF'' for Wirtinger flow, ``TWF'' for truncated Wirtinger flow, ``TAF'' for truncated amplitude flow, ``TRGrad'' for truncated Riemannian gradient descent. }\label{table:convergence}
\end{table}

\section{Main results}

We introduce tensors $\tT, \tS\in \left(\mathbb{R}^{2n}\right)^{\otimes 4}$ that will be useful to formulate our main results. Define $\tT$ as 
$$\tT:=\frac{1}{c} \sum_i(\va_i^+)^{\otimes 4},$$ 
where $c=m\sigma^4$ and $\sigma^2=\Var((\va_i^+)_1)$. We also define $\tS$ as 
$$\tS_{i_1,i_2,i_3,i_4}:=\tones_{i_1=i_2, i_3=i_4}+\tones_{i_1=i_3, i_2=i_4}+\tones_{i_1=i_4, i_2=i_3},$$
where $(\tones_{i_1=i_2, i_3=i_4})_{i_1,i_2,i_3,i_4}=1$ if $i_1=i_2$ and $i_3=i_4$ and $0$ else (and similarly for other tensors). \\  

The main results of this paper are given in \Cref{thm:1} and \Cref{thm:2}. \Cref{thm:1} gives a tensor based criterion for global benign landscape of $f$. 

\begin{theorem}\label{thm:1}
Assume $\|\tT-\tS\|_{\rm{op}}\le \delta_0$ for some constant $\delta_0>0$ small enough. The only local minima of $f$ defined in \cref{eq:realPR} are global minima $\vx^{\natural}e^{i\theta}$ and all saddle points of $f$ are strict\footnote{A strict saddle point is a critical point at which the Hessian has a strictly negative eigenvalue.}. 
\end{theorem}

The proof of \Cref{thm:1} is available in \Cref{subsec:landscape}. The main idea is to show that when $\|\tT-\tS\|_{\rm{op}}$ is small, the objective function $f$ is restricted strongly convex near the ground truth, has indefinite Hessian near the origin, and has nonzero gradient elsewhere. 

With \Cref{thm:1} at hand, we can derive our global convergence result in \Cref{thm:2}. 

\begin{theorem}\label{thm:2}
Assume $\|\tT-\tS\|_{\rm{op}}\le \delta_0$ for some constant $\delta_0>0$ small enough. For almost every initial point $\vx_0^+\in\mathbb{R}^{2n}$, there exists $\bar{\alpha}>0$ such that for any step sizes $(\alpha_k)_{k\in\mathbb{N}}$ satisfying $\sum_{k=0}^\infty\alpha_k=\infty$ and $\alpha_k\in(0,\bar{\alpha}],\,\forall k\in\mathbb{N}$, the gradient descent algorithm \cref{eq:gd_realPR} converges to a global minimizer of $f$ defined in \cref{eq:realPR}. 
\end{theorem}

The proof of \Cref{thm:2} is available in \Cref{subsec:convergence}. The main idea is to show that $f$ has bounded gradient trajectories, and then use the general convergence result \citep{josz2023global} together with \Cref{thm:1} to conclude the global convergence of gradient descent with random initialization for phase retrieval problem. 



The rest of this section will be organized as follows. In \Cref{subsec:equivalentform}, we introduce an equivalent formulation of $f$ as an inner product of some tensors. In \Cref{subsec:landscape}, we provide necessary ingredients for proving \Cref{thm:1}. In \Cref{subsec:convergence}, we show that $f$ has bounded gradient trajectories and deduce \Cref{thm:2}.

\subsection{Equivalent formulation of the problem}\label{subsec:equivalentform}

In this section, we provide an equivalent tensor based formulation for $f$. We write Hermitian product as a function of real vectors $\vx^+$ and $\vx^-$:
\begin{align}
    |\langle \va_i, \vx\rangle|^2=\langle \va_i^+, \vx^+\rangle^2+\langle \va_i^+, \vx^-\rangle^2. \label{eq:complexscalar-to-realscalar}
\end{align}

To make use of tensors, for $\vu,\vv\in \mathbb{R}^{2n}$, we develop the following product of scalar products as:
\begin{align}
    \langle \va_i^+,\vu\rangle^2\langle \va_i^+, \vv\rangle^2&=\sum_{i_1,i_2,i_3,i_4} (\va_i^+)_{i_1}(\va_i^+)_{i_2}(\va_i^+)_{i_3}(\va_i^+)_{i_4}\evu_{i_1}\evu_{i_2}\evv_{i_3}\evv_{i_4}\notag\\
    &=\langle (\va_i^+)^{\otimes 4}, \vu^{\otimes 2}\otimes \vv^{\otimes 2}\rangle. \label{eq:scalar-to-tensor}
\end{align}
Combining \Cref{eq:complexscalar-to-realscalar} and \Cref{eq:scalar-to-tensor}, we obtain an alternative expression of $f$ as in \Cref{prop:formulaf}. 

\begin{proposition}\label{prop:formulaf}
    $f$ can be written as 
    $$f(\vx^+)= \left\langle \sum_i(\va_i^+)^{\otimes 4}, \tU(\vx^+)\right\rangle = c\langle \tT,\tU(\vx^+)\rangle, $$
    where $\tU(\vx^+)$ is a tensor defined by 
    $$\tU(\vx^+):=\sum_k \varepsilon_k\vx_{1,k}^{\otimes 2}\otimes \vx_{2,k}^{\otimes 2}$$
    with $\vx_{1,k},\vx_{2,k}\in \{\vx^+, \vx^-, \vx^{\natural +}, \vx^{\natural -}\}$ and $\varepsilon_k=\pm 1 $. 
\end{proposition}

The proof of \Cref{prop:formulaf} is in \Cref{proof:equivForm}. We regard $\tU$ as a function of $\vx^+$ even if it involves $\vx^-$, because $\vx^-=\mM\vx^+$ is also a function of $\vx^+$. With the above alternative form of $f$, the assumption $\|\tT-\tS\|_{\rm{op}}\le \delta_0$ allows us to the landscape of $f$ by studying the landscape of its approximation $c\langle \tS, \tU(\vx^+)\rangle$. A detailed analysis of this approximation is given in \Cref{subsec:landscape}.

\subsection{Landscape results}\label{subsec:landscape}

In this subsection, we prove that all critical points of $f$ are either global minimizers or strict saddles. Given the fact that $\tT$ is close to $\tS$, we can expect that $f=c\langle \tT, \tU\rangle$ is close to $c \langle \tS, \tU\rangle$ in some sense. Therefore, our landscape analysis of $f$ is inspired by the landscape information of this approximation. \Cref{prop:expectation-f} gives an equivalent formula for this approximation function. 

\begin{proposition}\label{prop:expectation-f}
We have 
$$g(\vx^+):=c\langle \tS, \tU(\vx^+)\rangle=8c\left(\|\vx\|^4+\|\vx^{\natural}\|^4-|\langle \vx,\vx^{\natural}\rangle|^2-\|\vx\|^2\|\vx^{\natural}\|^2\right). $$
\end{proposition}

The proof of \Cref{prop:expectation-f} is in \Cref{proof:alternative-g}. Local minimizers of $g$ are global minimizers and other critical points of $g$ are strict saddles (see the analysis in \Cref{proof:unionisall}). Motivated by this, we prove that similar geometrical properties also hold for $f$. We define the following $3$ regions: 
\begin{align}&\mathcal{R}^1_{\delta_0}:=\left\{\vx~\bigg|~(\|\vx\|+\|\vx^{\natural}\|)^3\le \frac{\|\nabla g(\vx^+)\|}{C_1\delta_0 c}\right\}\label{eq:region1}, \\
&\mathcal{R}^2_{\delta_0}:=\left\{\vx~\bigg|~8\frac{|\langle \vx,\vx^{\natural}\rangle|^2}{\|\vx^{\natural}\|^4}+\left(4+\frac{1}{4}\delta_0C_2\right)\frac{\|\vx\|^2}{\|\vx^{\natural}\|^2} \le \left(4-\frac{1}{4}\delta_0 C_2\right)\right\}\label{eq:region-2}, \\
&\mathcal{R}^3_{\delta_0}:=\left\{\vx ~\bigg|~d(\vx, \mathcal{G})\le \frac{16-5C_2\delta_0}{192+4 C_2\delta_0} \right\}, 
\end{align}
where $c,\delta_0,C_1,C_2$ are constants and $\mathcal{G}$ is the set of global minimizers of $f$. 

Our goal is to show that 
\begin{itemize}
    \item $\mathcal{R}^1_{\delta_0}$ contains no critical points of $f$, 
    \item $\mathcal{R}^2_{\delta_0}$ only contains strict saddles of $f$, 
    \item $\mathcal{R}^{3}_{\delta_0}$ only contains global minimizers. 
\end{itemize}

\begin{figure}[H]
    \begin{center}
    \includegraphics[scale=0.8]{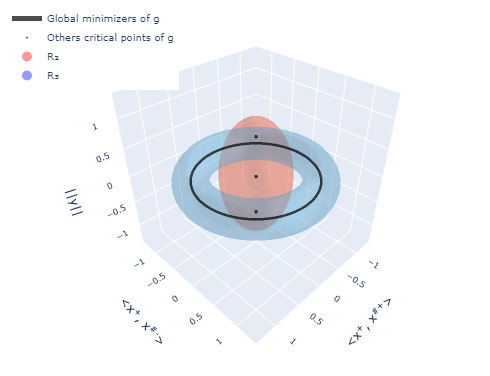}
    \caption{Overview of the regions considered in the proof}
    \label{fig:regions}
    \end{center}
\end{figure}

If we write $\vx^+=\vy+\mu\vx^{\natural +}+\nu \vx^{\natural -}$ with $\vy\in\mathbb{R}^{2n}$ orthogonal to $\vx^{\natural +},\vx^{\natural -} $, we can represent regions $\mathcal{R}^2_{\delta_0}, \mathcal{R}^3_{\delta_0}$ on a 3-d plot with coordinates  $(\|\vy\|, \langle \vx^+, \vx^{\natural +}\rangle,\langle \vx^+, \vx^{\natural -}\rangle) $ as shown in \Cref{fig:regions}. 

First we prove that $\mathcal{R}^1_{\delta_0}$ has no critical points. We bound the entries of the gradient of $\tU$ and use the concentration result $\|\tT-\tS\|_{\rm{op}}\le \delta_0$ to lower bound $\langle \nabla f(\vx^+), \nabla g(\vx^+)\rangle$ by $\|\nabla g(\vx^+)\|^2-C_1\|\nabla g(\vx^+)\|(\|\vx\|+\|\vx^{\natural}\|)^3$ and deduce \Cref{prop:R1-analysis}. 

\begin{proposition}\label{prop:R1-analysis}
    Assuming $\|\tT-\tS\|_{\rm{op}}<\delta_0$, we have for some absolute constant $C_1$, 
    $$ \|\nabla g(\vx^+)\|\ge C_1\delta_0 c(\|\vx\|+\|\vx^{\natural}\|)^3\implies \nabla f(\vx^+)\ne 0. $$
\end{proposition}

The proof of \Cref{prop:R1-analysis} is in \Cref{proof:R1-analysis}. This means that as long as the gradient of $g$ is large enough at $\vx^+$, then $\vx^+$ is not a critical point and of $f$ and thus $\mathcal{R}^1_{\delta_0}$ does not contain any critical point of $f$. 

To show that the critical points in $\mathcal{R}^{2}_{\delta_0}$ and $\mathcal{R}^{3}_{\delta_0}$ are either strict saddles or global minima respectively, we need to bound the distance between the Hessian of $f$ and $g$. We control the Hessian of the entries of $\tU$ and use the concentration $\|\tT-\tS\|_{\rm{op}}\le \delta_0$ to bound $\vu^T \mH_f \vu-\vu^T \mH_g \vu$ in \Cref{prop:hessian-distance}. 

\begin{proposition}\label{prop:hessian-distance}
Assume $\|\tT-\tS\|_{\rm{op}}<\delta_0$. For all vectors $u$ and for some absolute constant $C_2$, 
\begin{equation}
|\vu^T\mH_f(\vx^+)\vu-\vu^T\mH_g(\vx^+)\vu|< C_2\delta_0 c\|
\vu\|^2(\|\vx\|+\|\vx^{\natural}\|)^2. 
\label{eq:control-hessian}
\end{equation}
\end{proposition}

The proof of \Cref{prop:hessian-distance} is in \Cref{proof:hessian-dist}. To derive geometrical properties for $f$, we need an explicit formula for the Hessian of $g$, which can be easily deduced from \Cref{prop:expectation-f}. The result is given in \Cref{prop:hessian-g}. 

\begin{proposition}\label{prop:hessian-g}
    For all $\vx^+$ we have 
    \begin{equation}
    \mH_g(\vx^+) = 8c(8\vx^+(\vx^+)^{\rm T}+4\|\vx^+\|^2\mI_n - 2\vx^{\natural +}(\vx^{\natural +})^{\rm T}- 2\vx^{\natural -}(\vx^{\natural -})^{\rm T} - 2\|\vx^{\natural +}\|^2\mI_n  ). \label{eq:hessian}
\end{equation}
\end{proposition}

The proof of \Cref{prop:hessian-g} is in \Cref{proof:hessian-g}. Combining \Cref{prop:hessian-distance} and \Cref{prop:hessian-g}, if $\mH_g(\vx^+)$ has a sufficiently negative eigenvalue along some direction, then $\mH_f(\vx^+)$ will also have a negative eigenvlaue along the same direction. In region $\mathcal{R}^2_{\delta_0}$, we find that $\mH_g(\vx^+)$ indeed has a sufficiently negative eigenvalue along $\vx^{\natural +}$, which yields \Cref{prop:R2-analysis}. 
\begin{proposition}\label{prop:R2-analysis}
Assuming $\|\tT-\tS\|_{\rm{op}}<\delta_0$, we have 
$$8\frac{|\langle \vx,\vx^{\natural}\rangle|^2}{\|\vx^{\natural}\|^4}+\left(4+\frac{1}{4}\delta_0C_2\right)\frac{\|\vx\|^2}{\|\vx^{\natural}\|^2} \le \left(4-\frac{1}{4}\delta_0 C_2\right)\implies (\vx^{\natural +})^T\mH_f(\vx^+)\vx^{\natural +}<0. $$
\end{proposition}

The proof of \Cref{prop:R2-analysis} is in \Cref{proof:prop:R2-analysis}. It shows that all possible critical points in the region $\mathcal{R}^2_{\delta_0}$ are strict saddles. 

Finally we study the critical points in $\mathcal{R}^{3}_{\delta_0}$. We prove that all possible critical points in this region are global minimizers. To do so, we use the concept of restricted strong convexity similar to \cite{sun2018geometric}. The set of global minimizers is the circle $\mathcal{G}=\{\vx^{\natural}e^{i\theta}:\theta\in[0,2\pi)\}$, so we cannot expect $f$ to be strongly convex in this region. However, we can expect to have strong convexity in the directions orthogonal to $T_{\vx^+}(\mathcal{G})$ (the line tangent to the circle of global minimizers at $\vx^+$ where $\vx\in \mathcal{G}$). Again, we combine \Cref{prop:hessian-distance} and \Cref{prop:hessian-g} to prove that we indeed have restricted convexity and deduce that there are no critical points near $\mathcal{G}$ in \Cref{prop:R3-analysis}. 

\begin{proposition}\label{prop:R3-analysis}
Assume $\|\tT-\tS\|_{\rm{op}}<\delta_0$. If 
$$d(\vx, \mathcal{G})\le \frac{16-5C_2\delta_0}{192+4 C_2\delta_0} $$
then $\vx$ is a critical point if and only if $\vx$ is a global minimizer. 
\end{proposition}

The proof of \Cref{prop:R3-analysis} is in \Cref{proof:R3-analysis}. To conclude the landscape result, the last step is  prove that $\mathcal{R}^1_{\delta_0}\cup \mathcal{R}^2_{\delta_0}\cup \mathcal{R}^3_{\delta_0}=\mathbb{R}^{2n}$ for $\delta_0$ small enough. 

The key idea is that $\mathcal{R}^1_{\delta_0}$ is increasing in the sense of set inclusion as $\delta_0$ gets smaller. The only points in the complementary of $\mathcal{R}^1_{\delta_0}$ are those close to critical points. As $\mathcal{R}^2_{\delta_0}$ and $\mathcal{R}^3_{\delta_0}$ contain an open neighbourhood of the critical points, the complementary of $\mathcal{R}^1_{\delta_0}$ is strictly included in $\mathcal{R}^2_{\delta_0}\cup \mathcal{R}^{2}_{\delta_0}$ for small enough $\delta_0$. Therefore, we can finally conclude the landscape result for $f$ in \Cref{prop:union-regions}. 

\begin{proposition}\label{prop:union-regions}
For $\delta_0$ small enough, we have 
$$\mathcal{R}^1_{\delta_0}\cup \mathcal{R}^2_{\delta_0}\cup \mathcal{R}^3_{\delta_0}=\mathbb{R}^{2n}.$$
\end{proposition}

The proof of \Cref{prop:union-regions} is in \Cref{proof:unionisall}. Combining \Cref{prop:R1-analysis,prop:R2-analysis,prop:R3-analysis,prop:union-regions}, the result in \Cref{thm:1} follows immediately. 

\subsection{Convergence of gradient descent}\label{subsec:convergence}

In this subsection, we prove \Cref{thm:2} by using the fact that $f$ has bounded gradient trajectories and our landscape results in \Cref{thm:1}. 

We say that $f$ has bounded gradient trajectories if for every $\vx_0^+\in\mathbb{R}^{2n}$, the solution $\vx^+(\cdot)$ to the following initial value problem 
\begin{equation}\label{eq:GFpr}
    (\vx^+)'(t) = -\nabla f(\vx^+(t)), \quad \vx^+(0)=\vx_0^+
\end{equation}
satisfies that $\|\vx^+(t)\|\leq c_{\vx_0^+}$ for all $t\geq 0$, where $c_{\vx_0^+}$ is a constant dependent on $\vx_0^+$. 

The following \Cref{prop:bddflow} verifies that $f$ has bounded gradient trajectories for general positive semidefinite matrices $\{\mA_i\}_{i=1}^m$, not necessarily of the form $\va_i^+(\va_i^+)^{\rm T}+\va_i^-(\va_i^-)^{\rm T}$. 
\begin{proposition}\label{prop:bddflow}
    Let $\mA_i\in\mathbb{R}^{2n\times 2n}$ be symmetric positive semidefinite and $y_i\in\mathbb{R}$ for all $i=1,\ldots,m$. Then \cref{eq:realPR} has bounded subgradient trajectories. 
\end{proposition}

The proof of \Cref{prop:bddflow} is in \Cref{proof:bddflow}. Therefore, we are certified to apply results in \cite[Corollary 1]{josz2023global} and \Cref{thm:1} to conclude \Cref{thm:2}.

\section{Experiments}

In this section, we will show numerically the concentration of tensor $\tT$ and that $m=O(n)$ is enough to ensure $\|\tT-\tS\|_{\rm op}$ small . Then, we show the convergence of the loss function for gradient descent with a fixed $n$ versus different $m$ and that all trajectories do converge at a linear rate when $m$ is large enough even for large initializations. 

\begin{figure}[ht]
\begin{center}
    \includegraphics[scale=0.7]{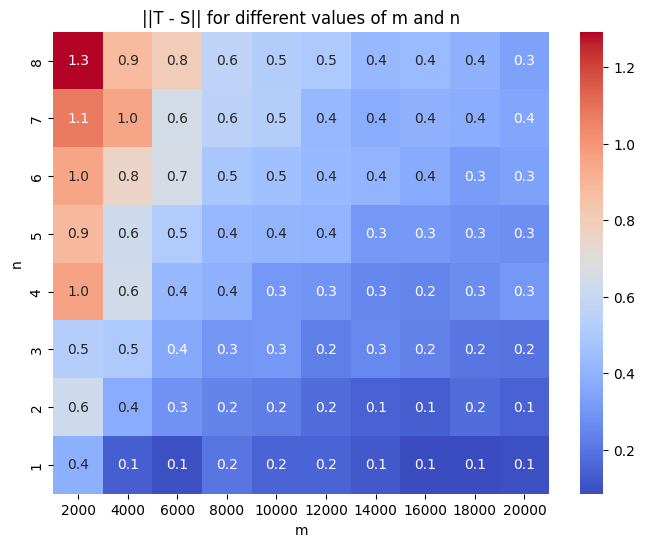}
    \caption{Values of $\|\tT-\tS\|_{\rm{op}}$ for different values of $m$ and $n$}
    \label{fig:operator-norms}
    \end{center}
\end{figure}

We generate a sample of $\ell=5$ sets of $m$ vectors $\{\va_i^{+}\}_{i=1}^{m}$ i.i.d and distributed as $\mathcal{N}(\vzero, \mI_{2n})$ for various $m,n$ and computed an approximation of $\sup_{\|\vu_1\|=\|\vu_2\|=\|\vu_3\|=\|\vu_4\|=1} \langle \tT, \vu_1\otimes \vu_2\otimes \vu_3\otimes \vu_4\rangle$ 
to estimate $\|\tT-\tS\|_{\rm{op}}$ and averaged the result over the $\ell$ samples. The result is displayed in \Cref{fig:operator-norms}. 

We can observe the concentration of $\tT$ around $\tS$ when $m$ is large enough for fixed $n$. Moreover, the boundary of the region $\|\tT-\tS\|_{\rm{op}}\le \delta_0$ for some $\delta_0$ is approximately linear: for $\delta_0=0.4$, the values of $m$ such that $\|\tT-\tS\|_{\rm{op}}\le \delta_0$ are approximately $m \ge 2000n$. This suggests that $m=O(n)$ samples are sufficient to ensure $\|\tT-\tS\|_{\rm{op}}$ is small enough. However, in practice, we don't need as much samples as $m=2000n$. As we will see in the next experiment, in practice the number of samples needed to recover the ground truth is about $m\approx 10n$ for almost all initial points. 

\begin{figure}[H]
    \centering
    \begin{minipage}{0.45\textwidth}
        \centering
        \includegraphics[width=\linewidth]{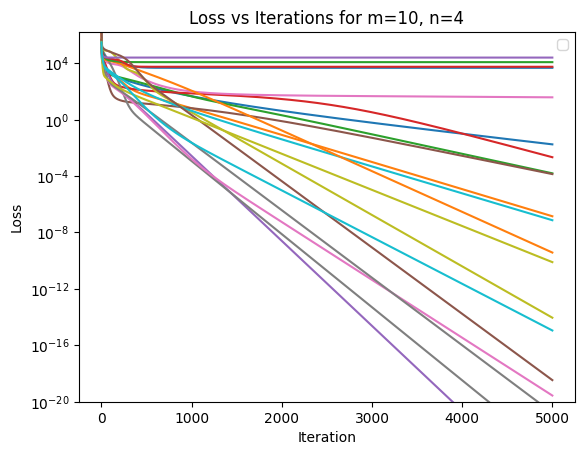}
       
    \end{minipage}
    \begin{minipage}{0.45\textwidth}
        \centering
        \includegraphics[width=\linewidth]{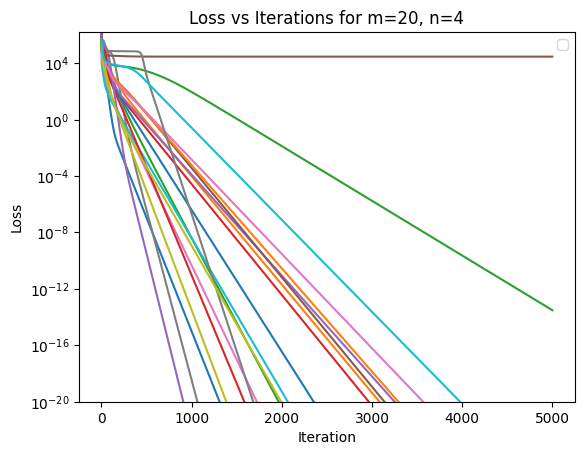}
        
    \end{minipage}

    \vspace{0.5cm}

    \begin{minipage}{0.45\textwidth}
        \centering
        \includegraphics[width=\linewidth]{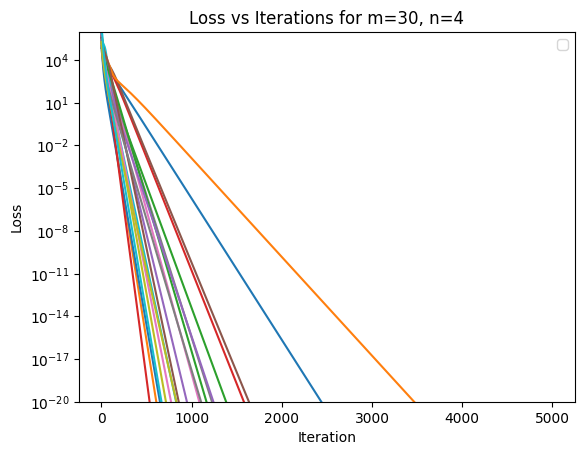}
        
    \end{minipage}
    \begin{minipage}{0.5\textwidth}
        \centering
        \includegraphics[width=\linewidth]{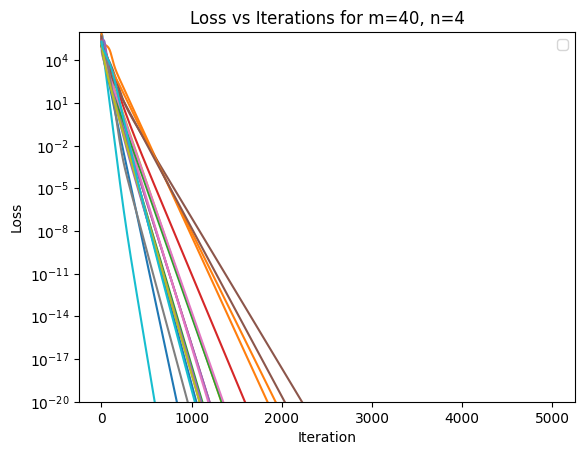}
    \end{minipage}

    \caption{Loss of $20$ trajectories for $n=4$ and $m\in [10,20,30,40]$}
    \label{fig:loss}
\end{figure}

In \Cref{fig:loss}, we plot the values $\tilde{f}=\frac{f}{m}$ at each iteration of vanilla gradient descent for $n=4$ and $m\in \{10,20,30,40\}$. We normalize the function $f$ by $\frac{1}{m}$ to have comparable magnitude for different values of $m$. For each value of $m$,  we chose $\{\va_i^+\}_{i=1}^{m}, \vx^{\natural +}$ i.i.d and distributed as $\mathcal{N}(\vzero, \mI_{2n})$. To study the influence of potentially large initializations, we select $\ell=20$ initial points following a uniform distribution on $[-10,10]^{n}$, a learning rate $\eta=5\cdot 10^{-5}$ and $N=5000$ iterations. 

As we may expect, some trajectories do not converge to a global minimum when $m$ is relatively small: for $m=10$ and $m=20$, there are $5$ and $1$ trajectories not converging to a global minimum respectively. However, when $m$ is large enough (for $m=30$ and $m=40$ in our experiment), all initializations converge to a global minimum of $\tilde{f}$ with a linear rate. 
Note that the rate of convergence also seems to be better when the number of samples increases.  

This suggests that the effective $K$ with which we have convergence of trajectories with high probability for $m\ge Kn$ should be around $K\approx 10$ as a rough estimate. Note that from \cite{fickus2014phase} theoretically there are at least $m=4n-2$ measurements needed to recover the ground truth vector, and using $m\approx 10n$ measurements is close to optimal. 

\section{Conclusion and discussion}

In this paper, we provided a tensor based criterion that guarantees global convergence of vanilla gradient descent with random initialization for phase retrieval problem. We first showed that the objective function has a benign global landscape if the criterion is verified. We also showed that the objective function has bounded gradient trajectories, which allows us to utilize the proposed landscape results and a general convergence result in literature to conclude our main convergence result. 

Finally, we discuss the limitation of our paper, which also leads to a potential future direction. From our proof, only a local linear convergence rate can be obtained once the iterates enter the region $\mathcal{R}_{\delta_0}^3$, but there is no information on how long the gradient descent algorithm will take to enter $\mathcal{R}_{\delta_0}^3$ starting from almost every point in $\mathcal{R}_{\delta_0}^1$ or $\mathcal{R}_{\delta_0}^2$. A more detailed analysis is needed for how long the iterates will stay in $\mathcal{R}_{\delta_0}^1$ and $\mathcal{R}_{\delta_0}^2$. For $\mathcal{R}_{\delta_0}^1$, the analysis is relatively simple. Since the gradient norm is lower bounded, the objective value will drop by a constant for every iteration. The key challenge is to ensure the iterates will not revisit $\mathcal{R}_{\delta_0}^2$ once leaving it. With the global landscape result provided in this paper, we expect to obtain a nearly linear iteration complexity result in a following paper.

\bibliography{references}
\bibliographystyle{abbrvnat}

\clearpage

\newpage
\appendix
\section{Appendix}

\subsection{Proof of equivalent formulation results}

\subsubsection{Proof  of \Cref{prop:formulaf}}\label{proof:equivForm}

\begin{proof} First we can express $|\langle \vx,\va_i\rangle|^2$ as a function of $\vx^+, \vx^-, \va_i^+$:
\begin{align}
    |\langle \vx,\va_i\rangle|^2 &=\left|(\mathrm{Re}(\vx)-i \mathrm{Im}(\vx))^T(\mathrm{Re}(\va_i)+i \mathrm{Im}(\va_i))\right|^2\notag\\
&=\left(\langle \mathrm{Re}(\vx), \mathrm{Re}(\va_i)\rangle + \langle \mathrm{Im}(\vx), \mathrm{Im}(\va_i)\rangle\right)^2+\left(\langle \mathrm{Re}(\vx), \mathrm{Im}(\va_i)\rangle -\langle \mathrm{Im}(\vx), \mathrm{Re}(\va_i)\rangle\right)^2\notag\\
&=\langle \vx^+,\va_i^+\rangle^2 +\langle \vx^-, \va_i^+\rangle^2.
\end{align}
We can then rewrite $f$ as:
\begin{align}
f(\vx)=\sum_i &\left(\langle \vx^+,\va_i^+\rangle^2 +\langle \vx^-, \va_i^+\rangle^2-\langle \vx^{\natural +},\va_i^+\rangle^2 -\langle \vx^{\natural -}, \va_i^+\rangle^2\right)^2\notag\\
=\sum_i&\bigg(\langle \vx^+, \va_i^+\rangle^4+\langle \vx^-, \va_i^+\rangle^4+\langle \vx^{\natural +}, \va_i^+\rangle^4+\langle \vx^{\natural -}, \va_i^+\rangle^4\notag\\
&+2\langle \vx^+,\va_i^+\rangle^2\langle \vx^-,\va_i^+\rangle^2+2\langle \vx^{\natural +},\va_i^+\rangle^2\langle \vx^{\natural -},\va_i^+\rangle^2\notag\\
&-2\langle \vx^+,\va_i^+\rangle^2\langle \vx^{\natural -},\va_i^+\rangle^2-2\langle \vx^-,\va_i^+\rangle^2\langle \vx^{\natural -},\va_i^+\rangle^2\notag\\
&-2\langle \vx^+,\va_i^+\rangle^2\langle \vx^{\natural +},\va_i^+\rangle^2-2\langle \vx^-,\va_i^+\rangle^2\langle \vx^{\natural +},\va_i^+\rangle^2\bigg).\label{eq:f-expanded}
\end{align}
From here we can notice that the scalar products can be written as scalar products with $(\va_i^+)^{\otimes 4}$, indeed for two vectors $\vu,\vv$, we have 
\begin{align}
\langle \vu,\va_i^+\rangle^2\langle \vv,\va_i^+\rangle^2&=\left(\sum_{i_1} \evu_{i_1}(\va_i^+)_{i_1}\right)\left(\sum_{i_2} \evu_{i_2}(\va_i^+)_{i_2}\right)\left(\sum_{i_3} \evu_{i_3}(\va_i^+)_{i_3}\right)\left(\sum_{i_4} \evu_{i_4}(\va_i^+)_{i_4}\right)\notag\\
&=\sum_{i_1,i_2,i_3,i_4} (\va_i^+)_{i_1}(\va_i^+)_{i_2}(\va_i^+)_{i_3}(\va_i^+)_{i_4}\evu_{i_1}\evu_{i_2}\evv_{i_3}\evv_{i_4}\notag\\
&=\langle (\va_i^+)^{\otimes 4}, \vu\otimes \vu\otimes \vv\otimes \vv \rangle. 
\end{align}

This allows to write $f$ in the following form:

\begin{align}
f(\vx) &=\sum_i \bigg\langle (\va_i^+)^{\otimes 4} , (\vx^+)^{\otimes 4}+(\vx^-)^{\otimes 4}+(\vx^{\natural +})^{\otimes 4}+(\vx^{\natural -})^{\otimes 4}\notag\\
&+2(\vx^+)^{\otimes 2}\otimes (\vx^-)^{\otimes 2}+2(\vx^{\natural +})^{\otimes 2}\otimes (\vx^{\natural -})^{\otimes 2}\notag\\
&-2(\vx^+)^{\otimes 2}\otimes (\vx^{\natural -})^{\otimes 2}-2(\vx^-)^{\otimes 2}\otimes (\vx^{\natural -})^{\otimes 2}\notag\\
&-2(\vx^+)^{\otimes 2}\otimes (\vx^{\natural +})^{\otimes 2}-2(\vx^-)^{\otimes 2}\otimes (\vx^{\natural +})^{\otimes 2} \bigg\rangle \notag\\
&=\bigg\langle \sum_i (\va_i^+)^{\otimes 4} , (\vx^+)^{\otimes 4}+(\vx^-)^{\otimes 4}+(\vx^{\natural +})^{\otimes 4}+(\vx^{\natural -})^{\otimes 4}\notag\\
&+2(\vx^+)^{\otimes 2}\otimes (\vx^-)^{\otimes 2}+2(\vx^{\natural +})^{\otimes 2}\otimes (\vx^{\natural -})^{\otimes 2}\notag\\
&-2(\vx^+)^{\otimes 2}\otimes (\vx^{\natural -})^{\otimes 2}-2(\vx^-)^{\otimes 2}\otimes (\vx^{\natural -})^{\otimes 2}\notag\\
&-2(\vx^+)^{\otimes 2}\otimes (\vx^{\natural +})^{\otimes 2}-2(\vx^-)^{\otimes 2}\otimes (\vx^{\natural +})^{\otimes 2} \bigg\rangle. 
 \end{align}

 We get the expected result by chosing:

\begin{align}
    \tU(\vx^+):&=(\vx^+)^{\otimes 4}+(\vx^-)^{\otimes 4}+(\vx^{\natural +})^{\otimes 4}+(\vx^{\natural -})^{\otimes 4}\notag\\
&+2(\vx^+)^{\otimes 2}\otimes (\vx^-)^{\otimes 2}+2(\vx^{\natural +})^{\otimes 2}\otimes (\vx^{\natural -})^{\otimes 2}\notag\\
&-2(\vx^+)^{\otimes 2}\otimes (\vx^{\natural -})^{\otimes 2}-2(\vx^-)^{\otimes 2}\otimes (\vx^{\natural -})^{\otimes 2}\notag\\
&-2(\vx^+)^{\otimes 2}\otimes (\vx^{\natural +})^{\otimes 2}-2(\vx^-)^{\otimes 2}\otimes (\vx^{\natural +})^{\otimes 2}.
 \end{align}
\end{proof}

\subsection{Proof of landscape results}

\subsubsection{Proof of \Cref{prop:expectation-f}}\label{proof:alternative-g}

\begin{proof}
In order to use tensor $\tS$, it is be useful to notice that:
\begin{align}
\langle \tS, \vu\otimes \vu\otimes \vv\otimes \vv\rangle &=\sum_{i_1,i_2,i_3,i_4}\tS_{i_1,i_2,i_3,i_4} \evu_{i_1}\evu_{i_2}\evv_{i_3}\evv_{i_4}\notag\\
&=\sum_{i_1=i_2\ne i_3=i_4} \evu_{i_1}\evu_{i_2}\evv_{i_3}\evv_{i_4}\notag\\
&+\sum_{i_1=i_3\ne i_2=i_4} \evu_{i_1}\evu_{i_2}\evv_{i_3}\evv_{i_4}\notag\\
&+\sum_{i_1=i_4\ne i_2=i_3} \evu_{i_1}\evu_{i_2}\evv_{i_3}\evv_{i_4}\notag\\
&+3\sum_{i_1=i_2=i_3=i_4} \evu_{i_1}\evu_{i_2}\evv_{i_3}\evv_{i_4}\notag\\
&=\|\vu\|^2\|\vv\|^2+2\langle \vu,\vv\rangle^2. \label{eq:formula-t}
\end{align}
 
Then we have:
\begin{align}
c\langle \tS, \tU(\vx^+)\rangle=&c \bigg(3\|\vx^+\|^4+3\|\vx^-\|^4+3\|\vx^{\natural +}\|^4+3\|\vx^{\natural -}\|^4\notag\\
&+2(2\langle \vx^+,\vx^-\rangle^2+\|\vx^+\|^2\|\vx^-\|^2)+2(2\langle \vx^{\natural +},\vx^{\natural -}\rangle^2 +\|\vx^{\natural +}\|^2\|\vx^{\natural -}\|^2)\notag\\
&-2(2\langle \vx^+, \vx^{\natural -}\rangle^2 +\|\vx^+\|^2\|\vx^{\natural -}\|^2)-2(2\langle \vx^-, \vx^{\natural -}\rangle^2 +\|\vx^-\|^2\|\vx^{\natural -}\|^2)\notag\\
&-2(2\langle \vx^-, \vx^{\natural +}\rangle^2 +\|\vx^-\|^2\|\vx^{\natural +}\|^2)-2(2\langle \vx^+, \vx^{\natural +}\rangle^2 +\|\vx^+\|^2\|\vx^{\natural +}\|^2)\bigg). 
\end{align}
We can notice that $\|\vx^+\|^2=\|\vx^-\|^2=\|\vx\|^2$, that $\langle \vx^+,\vx^-\rangle=0 $ (and same result holds for $\vx^{\natural}$) and that:
\begin{align}
&\langle \vx^+,\vx^{\natural +}\rangle ^2+\langle \vx^+,\vx^{\natural -}\rangle ^2+\langle \vx^-,\vx^{\natural +}\rangle ^2+\langle \vx^-,\vx^{\natural -}\rangle ^2\notag\\
=&2\left(\langle \mathrm{Re}(\vx),\mathrm{Re}(\vx^{\natural})\rangle +\langle \mathrm{Im}(\vx),\mathrm{Im}(\vx^{\natural})\rangle\right)^2+2\left(\langle \mathrm{Re}(\vx),\mathrm{Im}(\vx^{\natural})\rangle -\langle \mathrm{Im}(\vx),\mathrm{Re}(\vx^{\natural})\rangle\right)^2\notag\\
=&2|\langle \vx,\vx^{\natural}\rangle|^2. 
\end{align}
We can thus simplify the previous expression as:
$$c\langle \tS , \tU(\vx)\rangle=8c\left(\|\vx\|^4+\|\vx^{\natural}\|^4-|\langle \vx,\vx^{\natural}\rangle|^2-\|\vx\|^2\|\vx^{\natural}\|^2\right):=g(\vx^+). $$
\end{proof}

\subsubsection{Proof of \Cref{prop:R1-analysis}}\label{proof:R1-analysis}

\begin{proof}

To simplify notation, we denote $J=\{(i_1,i_2,i_3,i_4)~|~1\le i_1,i_2,i_3,i_4\le 2n\}$ and all summations over $j$ are for $j\in J$.  We also write:
$$\tU(\vx^+)=\sum_{k=1}^K \varepsilon_k \vx_{1,k}\otimes \vx_{2,k}\otimes \vx_{3,k}\otimes \vx_{4,k}$$
with $\vx_{i,k}\in \{\vx^+, \vx^-, \vx^{\natural +}, \vx^{\natural -}\}$ and $\varepsilon_k=\pm 1$. \\

Recall that:
\begin{align}
&f(\vx^+)=c\langle \tT, \tU(\vx^+)\rangle\notag\\
&g(\vx^+)=c\langle \tS, \tU(\vx^+)\rangle\notag\\
&\nabla f(\vx^+)=\sum_jc\tT_j\nabla \tU_j(\vx^+)\notag\\
&\nabla g(\vx^+)=\sum_jc\tS_j\nabla \tU_j(\vx^+). 
\end{align}
Using this we can write:
\begin{align}
    \langle \nabla f(\vx^+), \nabla g(\vx^+)\rangle&=\langle \nabla g(\vx^+), \nabla f(\vx^+)\rangle \notag\\
    &=\langle \nabla g(\vx^+), \sum_j c\tT_j \nabla \tU_j\rangle \notag\\
    &=\langle \nabla g(\vx^+), \sum_j c(\tT_j-\tS_j+\tS_j) \nabla \tU_j\rangle \notag\\
    &=\langle \nabla g(\vx^+), \nabla g(\vx^+)\rangle +\langle \nabla g(\vx^+), \sum_j c(\tT_j-\tS_j) \nabla \tU_j\rangle \notag\\
    &=\|\nabla g(\vx^+)\|^2 +\sum_i\sum_j (\nabla g(\vx^+))_i (\nabla \tU_j)_i c(\tT-\tS)_j \notag\\
    &=\|\nabla g(\vx^+)\|^2+ \sum_i\sum_k\sum_{i_1,i_2,i_3,i_4} (\nabla g(\vx^+))_ic(\tT-\tS)_{i_1,i_2,i_3,i_4} \varepsilon_k \notag\\
    &\times \bigg(\frac{\partial (\vx_{1,k})_{i_1}}{\partial \vx^+_{i}}(\vx_{2,k})_{i_2}(\vx_{3,k})_{i_3}(\vx_{4,k})_{i_4}\notag\\
    & +\frac{\partial (\vx_{2,k})_{i_2}}{\partial \vx^+_{i}}(\vx_{1,k})_{i_1}(\vx_{3,k})_{i_3}(\vx_{4,k})_{i_4}\notag\\
    &+\frac{\partial (\vx_{3,k})_{i_3}}{\partial \vx^+_{i}}(\vx_{1,k})_{i_1}(\vx_{2,k})_{i_2}(\vx_{4,k})_{i_4}\notag\\
    &+\frac{\partial (\vx_{4,k})_{i_4}}{\partial \vx^+_{i}}(\vx_{1,k})_{i_1}(\vx_{2,k})_{i_2}(\vx_{3,k})_{i_3}\bigg). 
\end{align}
If $\vx_{\ell,k}=\vx^+$, then $\frac{\partial (\vx_{\ell,k})_{i_{\ell}}}{\partial \evx_i}=\tones_{i=i_{\ell}}$. If $\vx_{\ell,k}=\vx^-$, then $\frac{\partial (\vx_{\ell,k})_{i_{\ell}}}{\partial \evx_i}=-\tones_{i=i_{\ell}+n}$ if $i_{\ell}\le n$ and $\tones_{i=i_{\ell}-n}$ otherwise. In each case we can write:
$$\sum_i (\nabla g(\vx^+))_i\frac{\partial (\vx_{\ell,k})_{i_{\ell}}}{\partial \evx_i}=(\mM_{\ell,k}\nabla g(\vx^+))_{i_{\ell}}$$
for some $\mM_{\ell,k}\in \left\{\mI_{2n}, \mO_{2n}, \begin{pmatrix}
    \mO_n & -\mI_n\\
    \mI_n & \mO_n
\end{pmatrix}\right\}$. Then we have 
\begin{align}
\langle \nabla f(\vx^+), \nabla g(\vx^+)\rangle&=\|\nabla g(\vx^+)\|^2+ \sum_k \bigg(c\langle \tT-\tS, \varepsilon_k \mM_{1,k}\nabla g(\vx^+)\otimes \vx_{2,k}\otimes \vx_{3,k}\otimes \vx_{4,k}\rangle \notag\\
&+c\langle \tT-\tS, \varepsilon_k  \vx_{1,k}\otimes\mM_{2,k}\nabla g(\vx^+)\otimes \vx_{3,k}\otimes \vx_{4,k}\rangle\notag\\
&+c\langle \tT-\tS, \varepsilon_k  \vx_{1,k}\otimes \vx_{2,k}\otimes\mM_{3,k}\nabla g(\vx^+)\otimes \vx_{4,k}\rangle\notag\\
&+c\langle \tT-\tS, \varepsilon_k  \vx_{1,k}\otimes \vx_{2,k}\otimes \vx_{3,k}\otimes \mM_{4,k}\nabla g(\vx^+)\rangle\bigg)\notag\\
&\ge  \|\nabla g(\vx^+)\|^2 - \sum_k \bigg(c\|\tT-\tS\|_{\rm{op}}\| \varepsilon_k \mM_{1,k}\nabla g(\vx^+)\|\|\vx_{2,k}\|\| \vx_{3,k}\|\vx_{4,k}|\|\notag\\
&+c\|\tT-\tS\|_{\rm{op}}\| \varepsilon_k \mM_{2,k}\nabla g(\vx^+)\|\|\vx_{1,k}\|\| \vx_{3,k}\|\vx_{4,k}|\|\notag\\
&+c\|\tT-\tS\|_{\rm{op}}\| \varepsilon_k \mM_{3,k}\nabla g(\vx^+)\|\|\vx_{1,k}\|\| \vx_{2,k}\|\vx_{4,k}|\|\notag\\
&+c\|\tT-\tS\|_{\rm{op}}\| \varepsilon_k \mM_{4,k}\nabla g(\vx^+)\|\|\vx_{1,k}\|\| \vx_{2,k}\|\vx_{3,k}|\|\bigg)\notag\\
&\ge  \|\nabla g(\vx^+)\|^2 - \sum_k 4\delta_0c (\|\vx\|+\|\vx^{\natural}\|)^3\|\nabla g(\vx^+)\|\label{eq:strict-gradient-inequality}\\
&\ge \|\nabla g(\vx^+)\|^2 - C_1\delta_0c (\|\vx\|+\|\vx^{\natural}\|)^3\|\nabla g(\vx^+)\|. 
\end{align}
Assume we have $ \|\nabla g(\vx^+)\|\ge C_1\delta_0 c(\|\vx\|+\|\vx^{\natural}\|)^3$, then in particular $\nabla g(\vx^+)\ne 0$, and $\vx^{\natural}\ne 0$, so inequality \cref{eq:strict-gradient-inequality} is strict and we must have $\langle \nabla f(\vx^+), \nabla g(\vx^+)\rangle >0$, and therefore $\nabla f(\vx^+)\ne 0$. 
    
\end{proof}

\subsubsection{Proof of \Cref{prop:hessian-distance}}\label{proof:hessian-dist}

We denote $J=\{(i_1,i_2,i_3,i_4)~|~1\le i_1,i_2,i_3,i_4\le 2n\}$ and all summations over $j$ are for $j\in J$.
\begin{proof}
    We have:
    \begin{align}
        &|\vu^T\mH_f(\vx^+)\vu-\vu^T\mH_g(\vx^+)\vu|\notag\\
        &=\left|\vu^T \left(c\sum_{j\in J} (\tT-\tS)_j \mH_{\tU_j}\right)\vu\right|\notag\\
        &=\left|\vu^T \left(c\sum_{j\in J}\sum_k \varepsilon_k (\tT-\tS)_j \mH_{(\vx_{1,k}\otimes \vx_{2,k}\otimes \vx_{3,k}\otimes \vx_{4,k})_j}\right)\vu\right|\notag\\
         &=\bigg|\vu^T \bigg(c\sum_{i_1,i_2,i_3,i_4}\sum_k \varepsilon_k (\tT-\tS)_{i_1,i_2,i_3,i_4} \notag\\
         &\sum_{\sigma \in \mathfrak{S}_4} \frac{1}{2}\frac{\partial (\vx_{\sigma(1),k})_{i_{\sigma(1)}}}{\vx^+_{i}}\frac{\partial (\vx_{\sigma(2),k})_{i_{\sigma(2)}}}{\vx^+_{i'}}( \vx_{\sigma(3),k})_{i_{\sigma(3)}}(\vx_{\sigma(4),k})_{i_{\sigma(4)}}\bigg)_{1\le i,i'\le 2n}\vu\bigg|\notag\\
         &=\bigg| c\sum_{i_1,i_2,i_3,i_4}\sum_k \varepsilon_k (\tT-\tS)_{i_1,i_2,i_3,i_4} \notag\\
         &\sum_{i,i'}\sum_{\sigma\in \mathfrak{S}_4} \frac{1}{2}\left(\evu_i\frac{\partial (\vx_{\sigma(1),k})_{i_{\sigma(1)}}}{\vx^+_{i}}\right)\left(\evu_{i'}\frac{\partial (\vx_{\sigma(2),k})_{i_{\sigma(2)}}}{\vx^+_{i'}}\right)( \vx_{\sigma(3),k})_{i_{\sigma(3)}}(\vx_{\sigma(4),k})_{i_{\sigma(4)}}\bigg|.\notag\\
    \end{align}
    As we saw previously in \Cref{proof:R1-analysis}, we can write for some $\mM_{\sigma(1),k}\in \left\{\mI_{2n}, \mO_{2n}, \begin{pmatrix}
    \mO_n & -\mI_n\\
    \mI_n & \mO_n
\end{pmatrix}\right\}$, 
    $$\sum_i\left(\evu_i\frac{\partial (\vx_{\sigma(1),k})_{i_{\sigma(1)}}}{\vx^+_{i}}\right)=(\mM_{\sigma(1),k}\vu)_{i_{\sigma(1)}}. $$
    Then we have:
    \begin{align}
        &|\vu^T\mH_f(\vx^+)\vu-\vu^T\mH_g(\vx^+)\vu|\notag\\
        &=\bigg| c\sum_{i_1,i_2,i_3,i_4}\sum_k \varepsilon_k (\tT-\tS)_{i_1,i_2,i_3,i_4} \notag\\
         &\frac{1}{2}\bigg((\mM_{1,k}\vu)_{i_1}(\mM_{2,k}\vu)_{i_2}(\vx_{3,k})_{i_3}(\vx_{4,k})_{i_4}+(\mM_{2,k}\vu)_{i_2}(\mM_{1,k}\vu)_{i_1}(\vx_{3,k})_{i_3}(\vx_{4,k})_{i_4}+\notag\\
         &\dots+ (\mM_{4,k}\vu)_{i_4}(\mM_{3,k}\vu)_{i_3}(\vx_{2,k})_{i_2}(\vx_{1,k})_{i_1}\bigg)\bigg|\notag\\
         &=\bigg|c\sum_k \varepsilon_k \bigg\langle \tT-\tS, \notag\\
         &\frac{1}{2}\bigg((\mM_{1,k}\vu)\otimes (\mM_{2,k}\vu)\otimes (\vx_{3,k})\otimes (\vx_{4,k})+(\mM_{1,k}\vu) \otimes(\mM_{2,k}\vu)\otimes (\vx_{3,k})\otimes (\vx_{4,k})+\notag\\
         &\dots+ (\vx_{1,k})\otimes (\vx_{2,k})\otimes (\mM_{3,k}\vu)\otimes (\mM_{4,k}\vu)\bigg)\bigg\rangle\bigg|\notag\\
          &\le  c\sum_k \|\tT-\tS\|_{\rm{op}} \sum_{\sigma\in \mathfrak{S}_4} \frac{1}{2}\|\mM_{\sigma(1),k}\vu\|\|\mM_{\sigma(2),k}\vu\| \|\vx_{\sigma(3),k}\|\|\vx_{\sigma(4),k}\|\notag\\ 
          &\le  c\sum_k  \|\tT-\tS\|_{\rm{op}} \sum_{\sigma\in \mathfrak{S}_4} \frac{1}{2}\|\vu\|^2(\|\vx\|+\|\vx^{\natural}\|)^2\notag\\ 
          &<  c\sum_k  \delta_0 12\|\vu\|^2(\|\vx\|+\|\vx^{\natural}\|)^2\notag\\
          &\le  C_2 \delta_0c\|\vu\|^2(\|\vx\|+\|\vx^{\natural}\|)^2. 
    \end{align}
\end{proof}

\subsubsection{Proof of \Cref{prop:hessian-g}}\label{proof:hessian-g}

\begin{proof}
We can notice that $\vx^-=\mM\vx^+$ where  
\begin{equation*}
    \mM := \begin{bmatrix}
        \mO_n & -\mI_n \\
        \mI_n & \mO_n
    \end{bmatrix}.
\end{equation*}
It is also well known that:
\begin{equation*}
    \mH_{\|\vx\|^2}(x) = 2\mI_n, \quad \mH_{\|\mM \vx\|^4}(\vx) = 8\vx\vx^{\rm T} + 4\|\vx\|^2\mI_n
\end{equation*}
\begin{equation*}
    \mH_{\langle \mM \vx,\va\rangle^2}(\vx) = 2(\mM \va)(\mM \va)^{\rm T}. 
\end{equation*}
Thus, the Hessian of $g$ with respect to $\vu^+$ is given by 
\begin{equation*}
    \mH_g(\vx^+) = 8c(8\vx^+(\vx^+)^{\rm T}+4\|\vx^+\|^2\mI_n - 2\vx^{\natural +}(\vx^{\natural +})^{\rm T}- 2\vx^{\natural -}(\vx^{\natural -})^{\rm T} - \|\vx^{\natural +}\|^2\mI_n - \|\vx^{\natural -}\|^2\mI_n ).
\end{equation*}
As $\|\vx^{\natural +}\| =\|\vx^{\natural -}\|$, we get the desired result. 
\end{proof}

\subsubsection{Proof of \Cref{prop:R2-analysis}}\label{proof:prop:R2-analysis}
\begin{proof}
Using \Cref{prop:hessian-g} and \Cref{prop:hessian-distance}, we can write:

\begin{align}
    (\vx^{\natural +})^T\mH_f(\vx^+)\vx^{\natural +}&< 8c\left(8\langle \vx^+,\vx^{\natural +}\rangle^2+4\|\vx^+\|^2\|\vx^{\natural +}\|^2-2\|\vx^{\natural +}\|^2-2\langle \vx^{\natural +}\vx^{\natural -}\rangle ^2-\|\vx^{\natural +}\|^4-\|\vx^{\natural +}\|^2\|\vx^{\natural -}\|^2\right)\notag\\
    &+C_2\delta_0c \|\vx^{\natural}\|^2(\|\vx\|+\|\vx^{\natural}\|)^2\notag\\
    &\le 8c\left(8\left(\langle \mathrm{Re}(\vx^{\natural}), \mathrm{Re}(\vx)\rangle+\langle \mathrm{Im}(\vx^{\natural}), \mathrm{Im}(\vx)\rangle\right)^2+4\|\vx\|^2\|\vx^{\natural}\|^2-2\|\vx^{\natural}\|^4 -\|\vx^{\natural}\|^4-\|\vx^{\natural}\|^4\right)\notag\\
     &+2C_2\delta_0c\|\vx^{\natural}\|^2 (\|\vx\|^2+\|\vx^{\natural}\|^2)\notag\\
    &\le 8c (8|\langle \vx,\vx^{\natural}\rangle|^2+4\|\vx\|^2\|\vx^{\natural}\|^2-4\|\vx^{\natural}\|^4)+2C_2\delta_0c \|\vx^{\natural}\|^2(\|\vx\|^2+\|\vx^{\natural}\|^2). 
\end{align}

Then as long as the following condition is satisfied, we have $(\vx^{\natural})^T\mH_f\vx^{\natural}<0$:
\begin{align}
    &8c (8|\langle \vx,\vx^{\natural}\rangle|^2+4\|\vx\|^2\|\vx^{\natural}\|^2-4\|\vx^{\natural}\|^4) +2\delta_0 C_2c(\|\vx\|^2\|\vx^{\natural}\|^2+\|\vx^{\natural}\|^4)\le 0\notag\\
    \iff & 8\frac{|\langle \vx,\vx^{\natural}\rangle|^2}{\|\vx^{\natural}\|^4}+\left(4+\frac{1}{4}\delta_0C_2\right)\frac{\|\vx\|^2}{\|\vx^{\natural}\|^2} \le \left(4-\frac{1}{4}\delta_0 C_2\right).\notag\\
\end{align}

\end{proof}

\subsubsection{Proof of \Cref{prop:R3-analysis}}\label{proof:R3-analysis}

\begin{proof}
We borrow the notations from \cite{sun2018geometric} and we let:
$$\phi(\vx)=\mathrm{argmin}_{\theta\in[0,2\pi)} \|\vx-\vx^{\natural}e^{i\theta}\|$$
(And for $\vx\in \mathrm{Span}\{\vx^{\natural}, i\vx^{\natural}\}^{\perp}$ we just define $\phi(\vx)=0$). Intuitively $\phi(\vx)$ is the phase of $\vx$ in the real $2d$ plane spanned by $\vx^{\natural}$. This means we can write in particular:
$$\vx^+=r\left(\cos(\phi(\vx))\vx^{\natural +} +\sin(\phi(\vx))\vx^{\natural -}\right)+\vv$$
with $\vv$ orthogonal to $\mathrm{Span}\{\vx^{\natural +},\vx^{\natural -}\}$. 

Now consider $\vu=r'\left(\cos(\phi(\vx))\vx^{\natural +}+\sin(\phi(\vx))\vx^{\natural -}\right)+\vv'$ such that $\|\vu\|=1$ and write $(\vx^{\natural})^{\phi(\vx)}=\cos(\phi(\vx))\vx^{\natural +}+\sin(\phi(\vx))\vx^{\natural -}$, we have 
\begin{align}
    \vu^T\mH_f(\vx^+)\vu&> \vu^T\mH_g(\vx^+)\vu -C_2\delta_0c(\|\vx^{\natural}\|+\|\vx\|)^2\qquad\qquad \text{(From \Cref{prop:hessian-distance})}\notag\\
    &=8c(8\langle \vu,\vx^+\rangle^2+4\|\vx^+\|^2\|\vu\|^2-2\langle \vx^{\natural +}, \vu\rangle ^2-2\langle \vx^{\natural -},\vu\rangle ^2-2\|\vx^{\natural +}\|^2\|\vu\|^2\notag\\
    &-C_2\delta_0c\|\vu\|^2(\|\vx^{\natural}\|+\|\vx\|)^2 \qquad\qquad\text{(From \Cref{prop:hessian-g})}\notag\\
    &=\|\vu\|^28c\bigg(\left(8\frac{\langle \vu,\vx^+\rangle^2}{\|\vu\|^2}+4\|\vx^+\|^2-2r'^2\|\vx^{\natural}\|^4-2\|\vx^{\natural}\|^2\right)-\frac{C_2}{8}\delta_0(\|\vx^{\natural}\|+\|\vx\|)^2\bigg)\notag\\
    &=  8c\bigg(\left(8\langle \vu,\vx^+\rangle^2+4\|\vx^+\|^2-2r'^2\|\vx^{\natural}\|^4-2\|\vx^{\natural}\|^2\right)-\frac{C_2}{8}\delta_0(\|\vx^{\natural}\|+\|\vx\|)^2\bigg)\notag\\
    &=  8c\bigg(\left(8\langle \vu,\vx^+-(\vx^{\natural})^{\phi(\vx)}+(\vx^{\natural})^{\phi(\vx)}\rangle^2+4\|\vx^+-(\vx^{\natural})^{\phi(\vx)}+(\vx^{\natural})^{\phi(\vx)}\|^2-2r'^2\|\vx^{\natural}\|^4-2\|\vx^{\natural}\|^2\right)\notag\\
    &-\frac{C_2}{8}\delta_0(\|\vx^{\natural}\|+\|\vx\|)^2\bigg)\notag\\
    &\ge   8c\bigg(\bigg(8\langle \vu,(\vx^{\natural})^{\phi(\vx)}\rangle^2+16\langle \vu,\vx^+-(\vx^{\natural})^{\phi(\vx)}\rangle\langle \vu,(\vx^{\natural})^{\phi(\vx)}\rangle+4\|(\vx^{\natural})^{\phi(\vx)}\|^2\notag\\
    &+8\langle \vx^+-(\vx^{\natural})^{\phi(\vx)},(\vx^{\natural})^{\phi(\vx)}\rangle  -2\|\vx^{\natural}\|^2-2r'^2\|\vx^{\natural}\|^4\bigg)-\frac{C_2}{8}\delta_0(\|\vx^{\natural}\|+\|\vx-\vx^{\natural}e^{i\phi(\vx)}\|+\|\vx^{\natural} e^{i\phi(\vx)}\|)^2\bigg)\notag\\
    &\ge 8c\bigg(8r'^2\|\vx^{\natural}\|^4-16\|\vx-\vx^{\natural}e^{i\phi(\vx)}\|\|\vx^{\natural}\|+4\|\vx^{\natural}\|^2-8\|\vx-e^{i\phi(\vx)}\|\|\vx^{\natural}\|-2r'^2\|\vx^{\natural}\|^4-2\|\vx^{\natural}\|^2\notag\\
    & -\frac{C_2}{8}\delta_0\left(4\|\vx^{\natural}\|^2+4\|\vx-\vx^{\natural}e^{i\phi(\vx)}\|\|\vx^{\natural}\|+\|\vx-\vx^{\natural}e^{i\phi(\vx)}\|^2\right)\bigg)\notag\\
    &= 8c\bigg(6r'^2\|\vx^{\natural}\|^4-24\|\vx-\vx^{\natural}e^{i\phi(\vx)}\|\|\vx^{\natural}\|+2\|\vx^{\natural}\|^2\notag\\
    & -\frac{C_2}{8}\delta_0\left(4\|\vx^{\natural}\|^2+4\|\vx-\vx^{\natural}e^{i\phi(\vx)}\|\|\vx^{\natural}\|+\|\vx-\vx^{\natural}e^{i\phi(\vx)}\|^2\right)\bigg)\notag\\
    &\ge 8c\bigg(2\|\vx^{\natural}\|^2-24\|\vx-\vx^{\natural}e^{i\phi(\vx)}\|\|\vx^{\natural}\| -\frac{C_2}{8}\delta_0(4\|\vx^{\natural}\|^2+4\|\vx-\vx^{\natural}e^{i\phi(\vx)}\|\|\vx^{\natural}\|+\|\vx-\vx^{\natural}e^{i\phi(\vx)}\|^2)\bigg). \notag
\end{align}
Assume $\|\vx-\vx^{\natural}e^{i\phi(\vx)}\|\le \|\vx^{\natural}\|$. If we have 
$$2\|\vx^{\natural}\|^2-24\|\vx-\vx^{\natural}e^{i\phi(\vx)}\|\|\vx^{\natural}\| \ge \frac{C_2}{8}\delta_0(5\|\vx^{\natural}\|^2+4\|\vx-\vx^{\natural}e^{i\phi(\vx)}\|\|\vx^{\natural}\|), $$
then we have $\vu^T\mH_f\vu>0$. 

Let $\mathcal{G}$ be the set of global minimizers of $f$, define:

\begin{align}
    \mathcal{R}^3_{\delta_0}&=\left\{2\|\vx^{\natural}\|^2-24\|\vx-\vx^{\natural}e^{i\phi(\vx)}\|\|\vx^{\natural}\| \ge \frac{C_2}{8}\delta_0(5\|\vx^{\natural}\|^2+4\|\vx-\vx^{\natural}e^{i\phi(\vx)}\|\|\vx^{\natural}\|)\right\}\notag\\
    &\cap \left\{\|\vx-\vx^{\natural}e^{i\phi(\vx)}\|\le \|\vx^{\natural}\|\right\}\notag\\
    &=\left\{\vx ~|~ d(\vx, \mathcal{G})\le \min\left(1,\frac{2-\frac{5}{8}C_2\delta_0}{24+\frac{4}{8}C_2\delta_0}\right) \right\}\notag\\
     &=\left\{\vx ~|~ d(\vx, \mathcal{G})\le \frac{16-5C_2\delta_0}{192+4 C_2\delta_0} \right\}. 
\end{align}
For $\delta_0$ small enough, note that $R^3_{\delta_0}$ strictly contains the set of global minimizers $\vx^{\natural}e^{i\theta}$. \\

Assume that there is a local minimum at $\vx^+\in \mathbb{R}^3_{\delta_0}$. Write:

$$\gamma(t)=t \vx^+ + (1-t)\cos(\phi(\vx))\vx^{\natural +}+(1-t)\sin(\phi(\vx))\vx^{\natural -}.$$
Consider $h(t)=f(\gamma(t))$, we have (with $\vx^+=\vv+r\cos(\phi(\vx))\vx^{\natural +}+r\sin(\phi(\vx))\vx^{\natural -}$):
\begin{align}
    h''(t)&=\left(\frac{d\gamma(t)}{dt}\right)^T \mH_f(\gamma(t))\left(\frac{d\gamma(t)}{dt}\right)+\nabla f(\gamma(t))\cdot \frac{d^2(\gamma(t))}{dt^2}\notag\\
    &=\left(\frac{d\gamma(t)}{dt}\right)^T \mH_f(\gamma(t))\left(\frac{d\gamma(t)}{dt}\right). 
\end{align}

Notice that $\frac{d\gamma(t)}{dt}$ can be written as $r'\left(\cos(\phi(\vx))\vx^{\natural +}+\sin(\phi(\vx))\vx^{\natural -}\right)+\vv'$, therefore $h''(t)>0$, so $h$ is strictly convex. As $h(0)$ is a global min of $h$, $h(1)$ cannot be a critical of $h$, and therefore $\vx$ cannot be a critical point of $f$ if $\vx\ne \vx^{\natural}e^{i\theta}$ for some $\theta$.
\end{proof}

\subsubsection{Proof of \Cref{prop:union-regions}}\label{proof:unionisall}
\begin{proof}
Recall that:
$$\nabla g(\vx^+)=8c\left(4\vx^+\|\vx^+\|^2-2\vx^+\|\vx^{\natural +}\|^2-2\vx^{\natural +}\langle \vx^+, \vx^{\natural +}\rangle  -2\vx^{\natural -}\langle \vx^{\natural -}, \vx^+\rangle\right). $$
We have: 
$$\|\nabla g(\vx^+)\|^2\underset{\|\vx\| \to +\infty}{\sim} (32 c)^2 \|\vx\|^6. $$
Therefore from \cref{eq:region1}, if we take $\delta_0$ such that:
$$(32c)^2> (C_1\delta_0c)^2\iff \delta_0< \frac{32}{C_1}.$$
Then $\mathcal{R}^1_{\delta_0}$ contains all $\vx$ for $\|\vx\|$ large enough. Then there exists some $K$ compact such that for all $\delta_0\le \frac{16}{C_1}$, we have $\vx\in \mathcal{R}^1_{\delta_0}$ for all $\vx\not \in K$. Let $E$ be an open set containing the critical points $\mathcal{C}$ of $g$. $K\backslash E$ is a compact also, therefore the function $\frac{\|\nabla g(\vx^+)\|}{C_1c(\|\vx\|+\|\vx^{\natural}\|)^3}$ attains its minimum on $K\backslash E$. Take $\delta_0$ strictly below this minimum, and from \cref{eq:region1}, $\mathcal{R}^1_{\delta_0}$  contains all $K\backslash E$. \\
Now let's consider a critical point of $g(\vx^+)$. We can see that:

Assume $\vx^+$ is a critical point. \\\\
$\bullet$ If $\|\vx^+\|^2=\frac{1}{2}\|\vx^{\natural +}\|^2$, then we must have $\langle \vx^+,\vx^{\natural +}\rangle=0$ and $\langle \vx^+, \vx^{\natural -}\rangle =0$ as $\vx^{\natural +}$ is ortogonal to $\vx^{\natural -}$, which means that $\|\vx\|^2=\frac{1}{2}\|\vx^{\natural}\|^2$ and that $|\langle \vx,\vx^{\natural}\rangle|^2=\langle \vx^+,\vx^{\natural +}\rangle^2+\langle \vx^+, \vx^{\natural -}\rangle^2 =0$. In that case using \Cref{prop:R2-analysis}, $\vx$ is strictly in $\mathcal{R}^2_{\delta_0}. $ for $\delta_0$ small enough. \\\\
$\bullet$ If $\|\vx^+\|^2\ne \frac{1}{2}\|\vx^{\natural +}\|^2$, then $\vx^+$ must be in $\mathrm{Span}\{\vx^{\natural +},\vx^{\natural -}\}$, lets write $\vx^+=\mu \vx^{\natural +}+\nu \vx^{\natural -}$. We must have:
$$\begin{cases}
    4\mu (\mu^2+\nu^2)-4\mu=0\\
    4\nu (\mu^2+\nu^2)-4\nu=0
\end{cases}$$
which gives either $\mu^2+\nu^2=1$, or $\mu=\nu=0$. If $\mu=\nu=0$, then $\vx=0$, which is strictly in $\mathcal{R}^2_{\delta_0}$ for $\delta_0$ small enough. If $\mu^2+\nu^2=1$, then $\vx=e^{i\theta }\vx^{\natural}$, which is strictly in $\mathcal{R}^3_{\delta_0}$ for $\delta_0$ small enough.\\\\
Therefore the set of critical points $\mathcal C$ is in $E:=\mathrm{Int}(\mathcal{R}^2_{\delta_0}\cup \mathcal{R}^3_{\delta_0}) $, the interior of $\mathcal{R}^2_{\delta_0}\cup \mathcal{R}^3_{\delta_0}$ for $\delta_0$ small enough, and from the previous argument we must have for $\delta_0$ small enough:
$$\mathcal{R}^1_{\delta_0}\cup \mathcal{R}^2_{\delta_0}\cup \mathcal{R}^3_{\delta_0}=\mathbb{R}^{2n}.$$

\end{proof}

\subsection{Proof of global convergence}

\subsubsection{Proof of \Cref{prop:bddflow}}\label{proof:bddflow}
\begin{proof}
Since $\mA_i$ is real symmetric and positive semidefinite for all $i=1,\ldots,m$, by orthogonal decomposition, we can write $\mA_i=\sum_{j=1}^r\lambda_{ij}\vv_{ij}\vv_{ij}^{\rm T}$, where $(\vv_{ij})_{j=1}^r$ are orthonormal, for some $1\le r\le n$ and $\lambda_{ij}>0$ for all $j=1,\ldots,r$. Define 
\begin{equation*}
    V:=\text{Span}\{\vv_{ij}:i=1,\ldots,m,\,j=1,\ldots,r\}. 
\end{equation*}
Notice that 
\begin{equation*}
    \nabla f(\vx^+)=4\sum_{i=1}^m\left(\langle \mA_i\vx^+,\vx^+\rangle-b_i\right)\mA_i\vx^+=4\sum_{i=1}^m\sum_{j=1}^r\lambda_{ij}\langle\vv_{ij},\vx^+\rangle\left(\langle \mA_i\vx^+,\vx^+\rangle-b_i\right)\vv_{ij} \in V.
\end{equation*}
Therefore, $\nabla f(\vx^+)\in V$ for all $\vx^+\in\mathbb{R}^{2n}$. Denote $V^{\perp}$ as the orthogonal complement of the subspace $V$, then for any given initial point $\vx_0^+$, the solution $\vx^+(\cdot)$ to \cref{eq:GFpr} can be decomposed as $\vx^+(t)=\vx^+_V(t)+\vx^+_{V^\perp}(t)$, where $\vx^+_V(t)\in V$ and $\vx^+_{V^\perp}(t)\in V^\perp$ for all $t\ge 0$. Note that $(\vx^+)'(t)=-\nabla f(\vx^+(t)) \in V$, thus $\vx^+_{V^\perp}(t)\equiv \vx^+_{V^\perp}(0)$ and we write $\vx^+(t)=\vx^+_V(t)+\vx^+_{V^\perp}(0)$ for all $t\ge 0$. Since $f(\vx^+)$ is a decreasing function over $t\ge 0$, 
\begin{align*}
    \sum_{j=1}^r\lambda_{ij}\langle \vv_{ij}, \vx^+(t)\rangle^2 &= |\langle \mA_i\vx^+(t),\vx^+(t)\rangle| \le \sqrt{2(\langle \mA_i\vx^+(t),\vx^+(t)\rangle-b_i)^2+2b_i^2} \\
    &\le \sqrt{2f(\vx^+(t))+2b_i^2} \le \sqrt{2f(\vx^+_0)+2b_i^2}. 
\end{align*}
Recall that $\lambda_{ij}>0$ for all $j=1,\ldots,r$, hence $\langle \vv_{ij},\vx^+(t)\rangle$ is bounded over $t\ge 0$ and so is $\langle \vv_{ij},\vx^+_V(t)\rangle$. As $\mathrm{Span}\{\vv_{ij}:i=1,\ldots,m,\,j=1,\ldots,r\}=V$, we can extract a basis of vectors $\vv_{ij}$ to form a basis of $V$, and denote this basis as $\{\vu_\ell:\ell=1,\ldots,d\}$. Then one can write $\vx^+_V(t)=\sum_{\ell=1}^d \zeta_\ell(t) \vu_\ell$. Notice that for each $\ell=1,\ldots,d$, there must exist $(i_\ell,j_\ell)$ such that $\vv_{i_\ell j_\ell}=\vu_\ell$. Thus, 
\begin{equation*}
    \|\vx^+_V(t)\|^2 = \sum_{\ell=1}^d\zeta_\ell(t)^2 = \sum_{\ell=1}^d\langle \vu_\ell,\vx^+_V(t)\rangle^2 = \sum_{\ell=1}^d\langle \vv_{i_\ell j_\ell},\vx^+_V(t)\rangle^2 
\end{equation*}
is bounded over $t\ge 0$. Finally, $\vx^+(t)=\vx^+_V(t)+\vx^+_{V^\perp}(0)$ is bounded over $t\ge 0$. 
\end{proof}

\end{document}